\documentclass[12pt]{article}
\usepackage{amsfonts,amssymb,amsbsy,amsmath,amsthm,enumerate}
\usepackage{color}

\newtheorem{theorem}{Theorem}[section]

\newtheorem{lemma}[theorem]{Lemma}
\newtheorem{follow}[theorem]{Corollary}
\newtheorem{assumption}[theorem]{Assumption}
\newtheorem{pr}[theorem]{Proposition}
\theoremstyle{definition}
\newtheorem{remark}[theorem]{Remark}

\newcommand{\bel}{\begin{equation} \label}
\newcommand{\ee}{\end{equation}}

\newcommand{\pd}{\partial}

\newcommand{\bx}{{\bf x}}
\newcommand{\E}{{\mathcal E}}

\newcommand{\rd}{{\mathbb R}^{2}}
\newcommand{\re}{{\mathbb R}}

\newcommand{\eps}{{\varepsilon}}

\def\beq{\begin{equation}}
\def\eeq{\end{equation}}
\newcommand{\bea}{\begin{eqnarray}}
\newcommand{\eea}{\end{eqnarray}}
\newcommand{\beas}{\begin{eqnarray*}}
\newcommand{\eeas}{\end{eqnarray*}}

{

\begin{document}
\begin{center}
{\LARGE \bf  Eigenvalue Asymptotics  in a Twisted  Waveguide}

\medskip

\today
\end{center}

\medskip

\begin{center}
{\sc Philippe Briet, Hynek Kova\v r\'{\i}k, Georgi  Raikov, Eric Soccorsi}\\
\end{center}

\begin{abstract}

We  consider a twisted  quantum wave guide i.e. a domain of the form
$\Omega_{\theta} : =  r_\theta \omega \times \re$ where $\omega$ is
a connected open and bounded subset of $\rd$ and $r_\theta=
r_\theta(x_3) $ is  a rotation by the angle $\theta(x_3)$ depending
on the longitudinal variable $x_3$. We are interested in the
spectral analysis of the Dirichlet Laplacian $H$ acting in
$\Omega_{\theta}$. We suppose that the derivative $\dot\theta$ of
the rotation angle can be written as $\dot\theta(x_3) =
\beta-\eps(x_3)$ with a positive constant $\beta$ and $\eps(x_3)
\sim L|x_3|^{-\alpha}$, $|x_3| \to \infty$. We show that if $L>0$
and  $\alpha \in (0,2)$, or if $L > L_0 > 0$ and $\alpha = 2$, then
there is an infinite sequence of discrete eigenvalues lying below
the infimum of the essential spectrum of $H$, and obtain the main
asymptotic term of this sequence.
\end{abstract}



%%%%%%%%%%%%%%%%%%%%%%%%%%%%%%%%%%%%%%%%%%%%%%%%%%%%%%%%%%%%%%%%%%%%%%%%%%
%%%%%%%%%%%%%%%%%%%%%%%%%%%%%%%%%%%%%%%%%%%%%%%%%%%%%%%%%%%%%%%%%%%%%%%%%%

\section{Introduction}
\label{intro} \setcounter{equation}{0}
 In this paper we deal with
the spectrum of the Dirichlet Laplacian in the twisted waveguide
$\Omega_{\theta} : =  r_\theta \omega \times \re$ where $\omega$
is a connected open and bounded subset of $\rd$ with sufficiently
regular boundary, and $r_\theta= r_\theta(x_3) $ is  a rotation by
the angle $\theta(x_3)$ depending on the longitudinal variable
$x_3$.

If the twisting of $\Omega_{\theta}$ is only local, then it does
not affect the essential spectrum of the Dirichlet Laplacian.
However, it does stabilise the discrete spectrum against possible
negative perturbations. Namely, it has been found out recently
that, provided $\omega$ is not rotationally symmetric, the local
twisting of $\Omega$ can be interpreted as a kind of a repulsive
perturbation, see \cite{BMT,EKK}. This has several consequences
such as the absence of weakly coupled bound states of
Schr\"odinger operators in twisted waveguides,
see \cite{EKK}, \cite{Kr}.

From this point of view  the effect of a local twisting of a
three-dimensional waveguide  is similar to the effect of a local
magnetic field in a two-dimensional waveguide, \cite{EK}.
Moreover, if the twisting is not local but constant, then it even
increases the infimum  of the essential spectrum of the Dirichlet
Laplacian, see \cite{ExK}. This is again reminiscent of the
situation in two-dimensional strips with constant magnetic field,
see \cite{BRS,GS}. On the other hand, any local decrease of a
constant twisting will induce at least one bound state of the
corresponding Dirichlet Laplacian, see \cite{ExK}.

 In the present paper we  study in more detail the properties of these
bound states induced by twisting. It is well known that the
Schr\"odinger operator with a slowly decaying potential possesses
infinitely many bound states,  and that the asymptotic distribution
of these bound states depends on the behaviour of the potential at
infinity, see e.g.~\cite{LS, RS4}. Our aim is to obtain analogous
asymptotic results for the bound states which are not induced by
an external potential, but by the twisting of the waveguide.

We start with the analysis of the unperturbed operator, which
corresponds to a constant twisting. This operator is
translationally invariant in the longitudinal direction and
therefore allows a fiber decomposition with  fiber operators which
have  purely discrete spectrum, see Subsection \ref{fiber}. As in
the model with constant magnetic field, \cite{BRS,GS}, we first
analyze the structure of the corresponding band functions. In
particular, we prove the existence of an effective mass   at the
bottom of the spectrum of the unperturbed operator, see Theorem
\ref{cor1}. Then we show that if the constant twisting is
perturbed by the function $\eps = \eps(x_3)$ which decays slowly
enough at infinity, then the resulting operator has infinitely
many discrete bound states accumulating from below at the infimum
of the essential spectrum. Finally, we employ the approach of
\cite{R} in order to study the rate of accumulation of these bound
states. Roughly speaking, our main result, Theorem \ref{good},
says that the rate of accumulation is determined by the rate of
decay of $\eps$ at infinity, and by the geometry of the
cross-section $\omega$.

%%%%%%%%%%%%%%%%%%%%%%%%%%%%%%%%%%%%%%%%%%%%%%%%%%%%%%%%%%%%%%%%%%%%%%%%%
%%%%%%%%%%%%%%%%%%%%%%%%%%%%%%%%%%%%%%%%%%%%%%%%%%%%%%%%%%%%%%%%%%%%%%%%%
\section{The Dirichlet Laplacian}
\label{prelim} \setcounter{equation}{0}
    \subsection{The general case}
  Let $ \omega$ be a  bounded  domain
in $ \rd$ with $C^2$-boundary. Furthermore, we suppose that
$\omega$ contains the origin of $\rd$.  Set $\Omega = \omega
\times \re$. For $\bx=(x_1,x_2,x_3)$ $\in
\Omega$ we write $\bx=(x_t,x_3)$ with $x_t=(x_1,x_2)$. Assume
that $ \theta \in C^1(\re)$ and its
derivative is bounded on $\re$. Define the twisted domain
$$
\Omega_\theta = \{ r_\theta(x_3) (\bx),\ \bx \in \Omega\}
$$
where
\[ r_\theta(x_3) = \left( \begin{array}{ccc}
\cos \theta(x_3) & \sin \theta(x_3) & 0 \\
-\sin \theta(x_3) & \cos \theta(x_3) & 0 \\
0 & 0& 1 \end{array} \right). \] We consider  the  Dirichlet
Laplacian  $-\Delta^D $ in ${\rm
  L}^2(  \Omega_\theta)$, i.e.~the self-adjoint operator  generated in  ${\rm L}^2(
\Omega_\theta)$ by the closed quadratic form
$$
\tilde{\mathcal Q}_\theta[f]= \int_{\Omega_\theta} \vert \nabla f
(\bx)\vert^2 d \bx,\quad f\in {\rm D}(\tilde{\mathcal Q}_\theta)=
{\rm H}_0^1 (\Omega_\theta).
$$
Here ${\rm H}_0^1(\Omega_\theta)$ denotes, as usual, the closure of
$C_0^{\infty}(\Omega_\theta)$ in the topology of the standard
first-order Sobolev space ${\rm H}^1(\Omega_\theta)$. Introduce the
transform
$$
(\mathcal{U} f)(\bx) = f \left( r_\theta(x_3) (\bx) \right), \quad \bx \in \Omega, \quad
f \in {\rm L}^2( \Omega_\theta).
$$
It is easy to see that $\mathcal{U}$ is  a unitary  operator from
${\rm L}^2( \Omega_\theta)$ onto ${\rm L}^2(\Omega)$.  Note also
that $ \mathcal{U}( {\rm H}_0^1 ( \Omega_\theta )) = {\rm  H}_0^1
(\Omega ).$ Set
 $$
 \nabla_{t}:  =( \partial_{1},
\partial_{2})^T, \quad \Delta_t : = \partial_{1}^2 +
\partial_{2}^2, \quad \partial_{\varphi} := x_1
\partial_2 - x_2 \partial_1.
$$
Denote by $\dot\theta$ the derivative of $\theta$. Define the
operator $H_{\dot{\theta}}$ as the self-adjoint operator generated
in ${\rm L}^2(\Omega)$ by the closed quadratic form
    \bel{q}
    {\mathcal Q}_{\dot{\theta}}[f] : = \tilde{\mathcal Q}_{\theta}[{\mathcal U}^{-1}f]
    = \int_{\Omega}
(\vert \nabla_{t} f \vert ^2 + \vert \dot{\theta}(x_3)
\partial_{\varphi} f  + \partial_3 f \vert^2 )\, d \bx, \quad
f \in {\rm  H}_0^1 (\Omega ) = : D(H_{\dot{\theta}}^{1/2}).
 \ee
 Evidently, $H_{\dot{\theta}} = \mathcal{U} (-\Delta^D)\,
 \mathcal{U}^{-1}$. By a straightforward computation we find
out that $H_{\dot{\theta}}$ acts on its domain as
    \bel{sept2}
H_{\dot{\theta}} = -\Delta_t -( \dot{\theta}(x_3)
\partial_{\varphi} + \partial_3)^2.
    \ee
Moreover, since  $\omega$ is a bounded domain, and we impose
Dirichlet boundary conditions, the operator $H_{\dot{\theta}}$ is
strictly positive, and hence boundedly invertible in ${\rm
L}^2(\Omega)$.
\begin{remark}  \label{disc}
If $\omega$ is a disk centered at the origin, then $\Omega_{\theta} = \Omega$ for any twisting $\theta$, and the operator
$H_{\dot\theta}$ is unitarily equivalent to $H_0$. Note that $H_0$ has purely absolutely continuous
spectrum (see   e.g.  \eqref{g50} below).
\end{remark}

%%%%%%%%%%%%%%%%%%%%%%%%%%%%%%%%%%%%%%%%%%%%%%%%%%%%%%%%%%%%%%%%%%%%%%%%%%
\subsection{Constant twisting}
\label{fiber}
    In this subsection we assume that the twisting is
constant,  i.e.  there exists $\beta \in \re$ such that $\dot{\theta}(x_3)
= \beta$ for every $x_3 \in \re$.\\
 Let ${\mathcal F}$ be the
partial Fourier transform with respect to $x_3$, i.e.
   % \bel{11}
   $$
({\mathcal F}u)(x_t,p): = \frac{1}{\sqrt{2\pi}} \int_{\re}e^{-i p
x_3} u(x_t,x_3)dx_3 , \quad (x_t,p) \in \omega \times \re.
    $$
    %\ee
    Due
to the  translational invariance in the $x_3$-direction, the
operator
${\hat H}_\beta : = {\mathcal F} H_\beta {\mathcal F}^*$
 admits a direct integral decomposition
 %\bel{10}
 $$
 {\hat H}_{\beta} : = \int_{\re}^\oplus h_\beta(p) dp,
 $$
 % \ee
 where
 \bel{sept2a}
h_\beta(p) : = -\Delta_t - (\beta
\partial_{\varphi} + ip)^2, \quad p \in \re,
\ee
 is the self-adjoint operator
generated in ${\rm L}^2(\omega)$ by the closed quadratic form
    $$
q_{p}[u] = q_{p,\beta}[u] :=\int_{\omega} \left( | \nabla_{t}
  u(x_t) |^2 + | (i
  \beta \partial_{\varphi} - p \right) u(x_t) |^2)dx_t \,,\quad
u \in {\rm H}_0^1(\omega).
$$
    Evidently, for each $p \in \re$ the quadratic form $q_p$
    induces an equivalent norm on ${\rm H}_0^1(\omega)$. Moreover,
    since the quadratic form
    ${\rm Re}\,i \int_{\omega} \pd_\varphi u \overline{u} dx_t$,
    $u \in {\rm H}_0^1(\omega),
    $ is relatively bounded with respect to the quadratic form
    $q_{0,\beta}[u]$ with zero relative bound, we find that $h_\beta(p)$, $p
    \in {\mathbb C}$, is a Kato analytic family of type B
(see \cite{K} or \cite[Theorem XII.10]{RS4}).\\
Further, since $\omega$ is a bounded domain with $C^2$-boundary,
the domain of the operator $h_\beta(p)$ is ${\rm H}^2(\omega) \cap {\rm
H}^1_0(\omega)$ for each $p \in \re$ (see e.g. \cite{A}).\\
Next, by the compactness of the embedding ${\rm H}^1_0(\omega)
\hookrightarrow {\rm L}^2(\omega)$, the spectrum of the operator
$h_\beta(p)$ is purely discrete. Let $\left
\{E_j(p)\right\}_{j=1}^{\infty} = \left
\{E_j(p, \beta)\right\}_{j=1}^{\infty}$ be the non-decreasing sequence of
the eigenvalues of $h_\beta(p)$, $p \in \re$.
Since $h_\beta(p)$ is a Kato analytic family, the functions $\re \ni p
\mapsto E_j(p) \in (0,\infty)$, $j \in {\mathbb N}$, are continuous piece-wise
analytic functions (see e.g. \cite[Theorem XII.12]{RS4}). The
mini-max principle easily yields
$$
E_j(p) = p^2 (1 + o(1)), \quad p \to \pm \infty.
$$
Therefore, the general theory of analytically fibred operators
(see e.g. \cite[Section XIII.16]{RS4}) implies that the spectrum
of the operator $H_\beta$ is purely absolutely continuous. In
summary, we have
    \bel{g50}
\sigma(H_\beta) = \sigma_{\rm ac}(H_\beta) = \cup_{j \in {\mathbb
N}} E_j(\re) = [\E, \infty),
    \ee
with
    \bel{g51}
    \E = \E(\beta) : = \min_{p \in \re} E_1(p, \beta).
    \ee
    Furthermore, for $p \in \re$ introduce a family
$\left\{\psi_j(x_t; p)\right\}_{j=1}^{\infty}$ of eigenfunctions
of the operator $h(p)$ which satisfy
    \bel{eideq} (h(p) \psi_j)(x_t;p) = E_j(p)\psi_j(x_t;p), \quad x_t \in \omega,
    \ee
and form an orthonormal basis in ${\rm L}^2(\omega)$. By the
embedding ${\rm H}^2(\omega) \hookrightarrow
C^l(\overline{\omega})$ with $l \in [0,1)$, we have
$\psi_j(\cdot;p) \in C^l(\overline{\omega})$, $j \in {\mathbb N}$;
in particular, the eigenfunctions $\psi_j(\cdot ;p)$ are
continuous and bounded on $\overline{\omega}$. Moreover, since the
operator $h_\beta(p)$ is strongly elliptic with coefficients in
 $C^\infty(\omega)$, we have $\psi_j(\cdot;p) \in C^\infty(\omega)$,  $j \in {\mathbb N}$
 (see e.g. \cite{A}, \cite{RS4}).

 Next, note that $h_\beta(0)$ is a strongly elliptic operator with
 real coefficients. Therefore, $E_1(0, \beta)$ is a simple eigenvalue of
 $h_\beta(0)$. Hence, we can choose $\psi_1(\cdot ; 0)$ to be real and
 positive in $\omega$ (see e.g. \cite[Theorem 8.38]{GT}). Moreover, there exists $\delta>0$
 such that the eigenvalue $E_1(p)$ is simple for $p \in [-\delta,\delta]$, and therefore the mapping $[-\delta,
 \delta] \ni p \mapsto E_1(p) \in (0,\infty)$ is analytic.
 Similarly,
 $\psi_1(\cdot;p)$ could be chosen so that the mapping $[-\delta,
 \delta] \ni p \mapsto \psi_1(\cdot ; p) \in {\rm H}^2(\omega)$ is
 analytic.

 At the end of this section, we establish for further references a criterion which guarantees that $\omega$ is a
 disk centered at the origin.

 \begin{pr} \label{septp1}
 Let $\omega \subset \rd$ be a  bounded domain with $C^2$-boundary $\partial \omega$.
 Assume that $\omega$ contains the origin.
 Then $\omega$ is a disk centered at the origin if and only if
    \bel{sept1}
    \|\partial_\varphi \psi_1(\cdot ; 0)\|_{{\rm L}^2(\omega)} = 0.
    \ee
 \end{pr}
 \begin{proof}
 For brevity set $\psi(x_t) = \psi_1(x_t;0)$, $x_t \in \omega$. Note that \eqref{sept1} implies $\pd_\varphi \psi = 0$
 identically in $\omega$ since $\psi \in C^{\infty}(\omega)$.\\
 (i) Assume that \eqref{sept1} holds true. Since $\psi$ is radial, it satisfies the equation
    \bel{sept3}
    (\Delta_t + \E) \psi = 0
    \ee
    in $\omega$ (see \eqref{sept2a}). Pick $\epsilon > 0$ such that the open
    disk $B_{\epsilon} : = \{x_t \in \rd \, | \, |x_t| < \epsilon\}$ is contained in $\omega$.
    Since $\psi$ is radial, regular, and satisfies \eqref{sept3} in $B_{\epsilon}$ we find that
    \bel{sept4}
    \psi(x_t) = \tilde{\psi}(x_t), \quad x_t \in B_{\varepsilon},
    \ee
    with
    $
    \tilde{\psi}(x_t) = c J_0(\E^{1/2} |x_t|)$, $x_t \in \rd$,
    where $J_0$ is the zeroth-order Bessel function (see e.g. \cite[Section 9.1]{AS}), and $c \neq 0$;
    if $c=0$, then the unique continuation principle (see e.g. \cite{JK}) would imply that $\psi = 0$ identically
    in $\omega$ which contradicts the fact that $\psi$ is an eigenfunction. Note that
    \bel{sept5}
    (\Delta_t + \E) \tilde{\psi} = 0
    \ee
    in $\rd$. Comparing \eqref{sept3} with \eqref{sept5}, and bearing in mind the unique continuation principle,
    we find that \eqref{sept4} holds for all $x_t \in \omega$. \\
    Let now $\{{\mathcal C}_{\alpha}\}$ be the set of the connected components of $\partial \omega$.
    Fix $\alpha$ and introduce the function $\varrho_\alpha$ by
    ${\mathcal C}_{\alpha} \ni x_t \mapsto \varrho_{\alpha}(x_t) : = |x_t| \in (0,\infty)$.
    Set ${\mathcal I}_\alpha : = \varrho_{\alpha}({\mathcal C}_\alpha)$. Since ${\mathcal C}_{\alpha}$ is connected,
    and $\varrho_{\alpha}$ is continuous, ${\mathcal I}_\alpha$ should be connected too, i.e. ${\mathcal I}_\alpha$ is
    a one-point set, or a bounded interval of positive length. Due to the Dirichlet boundary conditions,
    we  have $\psi(x_t) = 0$ for all $x_t \in {\mathcal C}_\alpha$, i.e.
    $J_0(\E^{1/2} r)=0$ for all $r \in {\mathcal I}_\alpha$. Since $J_0$ has at most a finite number of zeros
    on any bounded interval, all ${\mathcal I}_\alpha$ are one-point sets, i.e. all
    ${\mathcal C}_{\alpha}$ are arcs of circles centered at the origin. Since $\partial \omega \in C^2$,
     all ${\mathcal C}_\alpha$ are circles. Since  that $\omega$ is connected and contains the origin,
    it is a disk centered at the origin.\\
    (ii) Assume now that $\omega$ is a disk of radius $R \in (0,\infty)$, centered at the origin.
    Passing to polar coordinates $(r,\varphi)$, and writing
    $$
    u(r \cos{\varphi}, r \sin{\varphi}) = \frac{1}{\sqrt{2\pi}} \sum_{m \in {\mathbb Z}} e^{im\varphi} u_m(r), \quad u \in D(h_\beta(0)),
    $$
    we find that $h_\beta(0)$ is unitarily equivalent to
    $\oplus_{m \in {\mathbb Z}} {\mathcal H}_m$, where  ${\mathcal H}_m$, $m \in {\mathbb Z}$,
    is the operator generated in ${\rm L}^2((0,R);rdr)$ by the closure of the quadratic form
    $$\int_0^R \left(|u_m'(r)|^2 + m^2 (r^{-2} + \beta^2) |u_m(r)|^2\right) rdr, \quad u_m \in C_0^{\infty}(0,R).
    $$
    The spectra of all the operators ${\mathcal H}_m$ are discrete. Since
    $\E = \inf_{m \in {\mathbb Z}} \inf{\sigma({\mathcal H}_m)}$, and ${\mathcal H}_m \geq
    {\mathcal H}_0$ for all $m \in {\mathbb Z}$, we find that $\E$ coincides with
    the first eigenvalue of ${\mathcal H}_0$, which is simple. Let $\tilde{\psi}$
    be the real-valued eigenfunction of ${\mathcal H}_0$ which satisfies
    ${\mathcal H}_0 \tilde{\psi} = \E \tilde{\psi}$ and $\int_0^R \tilde{\psi}(r)^2 rdr = 1$.
    Since $\E$ is a simple eigenvalue of $h_\beta(0)$, it is not an eigenvalue of any ${\mathcal H}_m$
    with $m \neq  0$. Therefore, $\psi(x_t) = \tilde{\psi}(|x_t|)$, $x_t \in \omega$, i.e. $\psi$ is radial.
 \end{proof}

%%%%%%%%%%%%%%%%%%%%%%%%%%%%%%%%%%%%%%%%%%%%%%%%%%%%%%%%%%%%%%%%%%%%%%%%%%%%%%%
%%%%%%%%%%%%%%%%%%%%%%%%%%%%%%%%%%%%%%%%%%%%%%%%%%%%%%%%%%%%%%%%%%%%%%%%%%%%%%%
\section{Existence of an Effective Mass at the Origin}
 \label{effmass}   \setcounter{equation}{0}
 In this section  assume that $\dot{\theta}$ is equal to the constant $\beta \in \re$,
 and investigate the properties of the first band function $E_1(p)$, $p \in \re$.  We
 show  that  $E_1(0)$ attains its minimum value $\E$ at $p=0$
and only at $p=0$, and the minimum is non-degenerate, i.e. $E_1''(0)
> 0$. In this case, following a terminology established in the
solid-state physics (see e.g. \cite{KS}), we say that there
exists {\em an effective mass} at the origin. The existence of an effective mass
at the bottom of the absolutely continuous spectrum is an important
problem in the spectral analysis of many operators of the
mathematical physics (see e.g.
 \cite{KS, R, BSu, AP}).

For $\beta \in \re$ define
    \bel{epsilon}
\epsilon_{\omega}(\beta):= \frac{\beta^2 C_{\omega}}{1 +\beta^2
C_{\omega}}
    \ee
    where
    \bel{C}
    {C}_\omega : = \sup_{x_t \in \omega}(x_1^2 +
    x_2^2).
    \ee
    \begin{theorem}
    \label{cor1}
Let $\beta \in \re$.  Then \bel {masseff}
E_1(0,\beta) +
(1-\epsilon_{\omega}(\beta))\, p^2  \, \leq \, E_1(p,\beta)\, \leq E_1(0,\beta) + p^2, \quad  \, p \in \re.
    \ee
    \end{theorem}
To prove  the Theorem \ref {cor1} we  need  the following
technical result:
\begin{lemma}
\label{lm3} For all $p, \beta \in \re$ we have \bel{g61}
E_1(p,\beta) = E_1(0,\beta) + \inf_{0 \neq u \in
C_0^{\infty}(\omega)} \frac{\int_{\omega}\;\psi_1^2 ( | \nabla_{t}
u |^2 + | i \beta
  \partial_{\varphi} u - p u |^2 )dx_t }{\int_{\omega}\;\psi_1^2  \;| u |^2 dx_t },
  \ee
where $\psi_1 = \psi_1(\cdot;0)$.
\end{lemma}
\begin{proof}
The starting point of the proof is  the mini-max principle
    \bel{mima1}
E_1(p,\beta) =  \inf_{0 \neq v \in C_0^{\infty}(\omega)}
\frac{\int_{\omega} ( | \nabla_{t} v |^2 + | i \beta
  \partial_{\varphi} v - p
  v|^2 )\; dx_t}{\int_{\omega} | v |^2\; dx_t}.
    \ee
 Change the functional variable $v = {\psi_1} u$.  Note that $ \psi_1^{-1}
C_0^{\infty}(\omega)= C_0^{\infty}(\omega)$.
Integrating by parts, and bearing in mind that the function $\psi_1(\cdot;0)$ is real-valued, we easily get
$$ \int_{\omega} ( | \nabla_{t} v |^2 + | i \beta
\partial_{\varphi} v - p v |^2 ) dx_t =
$$
    \bel{g62}
\int_{\omega} \psi_1^2 ( | \nabla_{t} u |^2 + | i \beta
\partial_{\varphi} u
- p u |^2 ) dx_t + \int_{\omega} \psi_1 ( -\Delta_{t} \psi_1 -
\beta^2 \partial_{\varphi}^2 \psi_1 ) |u|^2 dx_t.
    \ee
Taking into account the eigenvalue equation
$$
-\Delta_{t} \psi_1 - \beta^2 \partial_{\varphi}^2 \psi_1 =
E_1(0,\beta) \psi_1,
$$
we find that (\ref{mima1}) and \eqref{g62} imply \eqref{g61}.
\end{proof}
 \begin{proof}[Proof of  Theorem \ref{cor1}]
We first prove the  lower bound. For $\beta =0$  it follows from
the Lemma \ref{lm3} and
$$
\int_{\omega} \psi_1^2 ( | \nabla_{t} u |^2 +p^2 | u |^2 ) dx_t
\geq p^2 \int_{\omega} \psi_1^2 | u|^2  dx_t, \quad u \in
C_0^{\infty}(\omega).
$$
Let us now consider the case $\beta \neq 0$. Pick $\eta \in (0,1)$.
Then we have
    $$ \int_{\omega} \psi_1^2 (|
\nabla_{t} u |^2  + | i \beta
\partial_{\varphi} u -
p u |^2) dx_t \geq
$$
$$
\int_\omega \psi_1^2 (| \nabla_{t} u |^2 +
\beta^2(1-\eta^{-1}) |\partial_\varphi u|^2 + (1-\eta) p^2 |u
|^2) dx_t \geq
$$
    \bel{mima5}
    (1+ \beta^2(1-\eta^{-1}) C_{\omega}) \int_{\omega}
\psi_1^2 | \nabla_{t} u |^2 dx_t + (1-\eta) p^2 \int_{\omega}
\psi_1^2 |u |^2 dx_t, \quad u \in C_0^{\infty}(\omega),
    \ee
    the constant $C_{\omega}$ being defined in \eqref{C}. Choose $\eta=
\epsilon_\omega \in (0,1)$, the constant $\epsilon_\omega$ being defined in \eqref{epsilon}.
     Then \eqref{mima5} implies
$$ \int_{\omega} \psi_1^2 (| \nabla_{t} u |^2 + | i \beta \partial_{\varphi}
u - p u |^2) dx_t \geq (1-\epsilon_\omega) p^2 \int_{\omega}
\psi_1^2 |u |^2 dx_t, $$ and  the lower bound in \eqref{masseff}
follows from Lemma \ref{lm3}.

 By the mini-max principle \eqref{mima1}, and the fact that $\psi_1$ is real-valued, we get
$$ E_1(p,\beta) \leq  \int_{\omega} ( | \nabla_{t} \psi_1 |^2
+ |i \beta
  \partial_{\varphi} \psi_1   - p
  \psi_1 |^2 )\; dx_t = E_1(0,\beta) + p^2,$$ which gives the
upper bound in \eqref{masseff}.
\end{proof}
Since $\epsilon_{\omega}(\beta) < 1$, and $E_1(p,\beta)$ depends
analytically on $p$ near the origin, Theorem \ref{cor1} entails the
following
\begin{follow} \label{gf1}
Let $\beta \in \re$. Then
we have $\pd_p E_1(0,\beta) = 0$ and
 \bel{nd1}
 \mu = \mu(\beta) : = \frac{1}{2} \pd_p^2 E_1(0,\beta) > 0,
 \ee
i.e. at the origin there exists an effective mass. As a corollary,
    \bel{nd2}
E_1(p,\beta) = E_1(0,\beta) + \mu(\beta)p^2 + O(p^3), \quad p \to 0.
    \ee
Moreover, for each $p \in \re$, $p \neq 0$, and $\beta \in \re$ we have
$E_1(p,\beta) > E_1(0,\beta)$.
\end{follow}

%%%%%%%%%%%%%%%%%%%%%%%%%%%%%%%%%%%%%%%%%%%%%%%%%%%%%%%%%%%%%%%%%%%%%%%%%%
%%%%%%%%%%%%%%%%%%%%%%%%%%%%%%%%%%%%%%%%%%%%%%%%%%%%%%%%%%%%%%%%%%%%%%%%%%
\section{Eigenvalue Asymptotics for the Dirichlet Laplacian with Non-Constant Twisting}
\label{s: main} \setcounter{equation}{0}
%%%%%%%%%%%%%%%%%%%%%%%%%%%%%%%%%%%%%%%%%%%%%%%%%%%%%%%%%%%%%%%%%%%%%%%%%%%%%%
\subsection{Perturbation of the twisting}
In this section the constant parameter $\beta > 0$ is fixed, and
we now consider the case where the derivative of the twisting is
given by
\begin{equation} \label{twist}
\dot\theta(x_3) = \beta- \eps(x_3)\,,
\end{equation}
where the perturbation $\eps\in {\rm L}^{\infty}(\re)$ satisfies
$\eps(x) \geq 0$, $x \in \re$, and $\lim_{x \to \pm \infty}
\eps(x) = 0$. By \eqref{sept2}, the operator $H_{\beta - \eps}$  can be
written as
    \bel{heps}
    H_{\beta - \eps} = H_{\beta} + W_{\eps, \beta}
    \ee
where
    \bel{g0} W_{\eps, \beta} : =  2\beta\eps\,
\pd_\varphi^2 + \pd_\varphi\, \eps\, \pd_3 +   \pd_3\,\eps
\pd_\varphi -\eps^2 \pd_\varphi^2 = - 2\beta G_0^* G_0 -  G_0^*
G_1 - G_1^* G_0 + G_2^* G_2,
    \ee
    and
    $$
    G_0 : = \eps^{1/2} \pd_\varphi, \quad G_1 : = \eps^{1/2}
    \pd_3, \quad G_2 : = \eps \pd_\varphi,
    $$
    are the operators defined  on $D(H_\beta^{1/2}) = {\rm H}_0^1(\Omega)$. Evidently, the operators
    $G_j H_\beta^{-1/2}$, $j=0,1,2$, are bounded
    in ${\rm L}^2(\Omega)$. Due to the Sobolev embedding theorems
    and the fact that $\eps(x) \to 0$ as $|x| \to \infty$, it
    is easy to see that
    the operators $G_j H_\beta^{-1}$, $j=0,1,2$, are compact
    in ${\rm L}^2(\Omega)$. Therefore, the operator
    $H_{\beta}^{-1/2} W_{\eps, \beta} H_{\beta}^{-1}$ is compact.
    Since the operators $H_{\beta - \eps}^{-1/2} H_{\beta}^{1/2}$ and,
    hence, $H_{\beta - \eps}^{-1} H_{\beta}^{1/2}$ are bounded, the resolvent difference
    $H_{\beta - \eps}^{-1} - H_\beta^{-1}$ is compact. By \cite[Theorem 4, Section 1, Chapter 9]{BSo2}, we have
    \bel{g25}
    \sigma_{\rm ess}(H_{\beta - \eps}) = \sigma_{\rm ess}(H_\beta).
    \ee
    \begin{remark} In a completely different context, a second-order perturbation with decaying
    coefficients which preserves the essential spectrum, has been
    considered in \cite{AADH}.
    \end{remark}
    Putting together \eqref{g25} and \eqref{g50}, we find that
    \bel{g60}
\sigma_{\rm ess}(H_{\beta - \eps}) = [\E, \infty)
    \ee
    where $\E$ is defined in \eqref{g51}.
%%%%%%%%%%%%%%%%%%%%%%%%%%%%%%%%%%%%%%%%%%%%%%%%%%%%%%%%%%%%%%%%%%%%%%
\subsection{Main result on the asymptotic distribution of the discrete spectrum of $H_{\beta - \eps}$}
It was shown in \cite{ExK} that a local decrease of the twisting
induces the existence of at least one bound state below $\E$
provided that $\omega$ is not radially symmetric. In other words,
even if $\eps$ is compactly supported then $H_{\beta - \eps}$ has
at least one discrete eigenvalue. Here we show that there is
actually an infinite number of bound states appearing below  $\E$
if the perturbation $\eps$ decays sufficiently slowly at infinity,
and describe the asymptotics  as $E \uparrow \E$ of the number of
bound states less than $E$.

We will impose on $\eps$ one of the following assumptions:

\begin{assumption} \label{estimate}
Suppose that $\eps \in C^1(\re)$, and there exist constants
$\alpha
> 0$ and $C>0$ such that
$$
0 \leq \eps(x) \leq C(1 + |x|)^{-\alpha},
$$
$$
|\dot{\eps}(x)| \leq C(1 + |x|)^{-\alpha-1}, \quad x \in \re.
$$
\end{assumption}

\begin{assumption} \label{asymptotics}
In addition to Assumption \ref{estimate}, suppose that
there exists a constant $L > 0$ such that
$$
\lim_{|x| \to \infty} |x|^{\alpha} \eps(x) = L.
$$
\end{assumption}

In order to formulate our main result, we need the following
notations. Let $T$ be a linear self-adjoint operator acting in a
given Hilbert space and such that $\tau:= \inf
\sigma_{ess}(T)>-\infty$. Then $N(T;t)$, $t \in (-\infty,\tau)$,
denotes the number of eigenvalues of $T$ lying on the interval
$(-\infty,t)$, and counted with the multiplicities.

\begin{theorem} \label{good}
(i) Suppose that Assumption \ref{asymptotics} holds  with
$\alpha \in (0,2)$. Then we have
    \bel{g26}
\lim_{\lambda\downarrow 0}\,  \lambda^{\frac 1\alpha-\frac 12}\,
N(H_{\beta-\eps}; \E -\lambda) = \frac{2}{\pi\alpha \sqrt{\mu}} \
\left(2\beta\, L\, \|\pd_\varphi\psi_1(\cdot, 0)\|^2_{{\rm
L}^2(\omega)} \right)^{\frac 1\alpha}\, B\left(\frac 32,\, \frac
1\alpha-\frac 12\right) \ee
where $\mu$ is defined in \eqref{nd1}, and $B$ is the Euler beta function.\\
(ii) Suppose that Assumption \ref{asymptotics} holds  with
$\alpha = 2$. Then we have
    \bel{g27}
\lim_{\lambda\downarrow 0}\,  |\ln{\lambda}|^{-1}\,
N(H_{\beta-\eps}; \E -\lambda) = \frac{1}{\pi}
\left(\frac{2\beta\, L}{\mu}\, \|\pd_\varphi\psi_1(\cdot,
0)\|^2_{{\rm L}^2(\omega)} - \frac{1}{4}\right)_+^{1/2}.
    \ee
If, moreover, $2\beta\, L\, \|\pd_\varphi\psi_1(\cdot,
0)\|^2_{{\rm L}^2(\omega)} < \frac{\mu}{4}$, then
    \bel{g1}
N(H_{\beta-\eps}; \E -\lambda) = O(1), \quad \lambda \downarrow 0.
    \ee
    (iii) Suppose that Assumption \ref{estimate} holds  with
$\alpha > 2$. Then we have \eqref{g1}.
\end{theorem}

\begin{remark}
Similarly to the potential perturbation case, the rate of
divergence of $N(H_{\beta-\eps}; \E-\lambda)$ when $\lambda
\downarrow 0$ is determined by the decay rate $\alpha$ of the
perturbation.  The asymptotic coefficients in \eqref{g26} and
\eqref{g27}  also depend on the limit $L$, the constant $\beta$,
and the geometry of the cross-section $\omega$ through the factor
$\|\pd_\varphi\psi_1(\cdot, 0)\|_{{\rm L}^2(\omega)}$. Note that
if $\|\pd_\varphi\psi_1(\cdot, 0)\|_{{\rm L}^2(\omega)} = 0$, then
the asymptotic coefficients in \eqref{g26} and \eqref{g27} vanish.
However, in this case Proposition \ref{septp1} implies that
$\omega$ is a disk centered at the origin. By Remark \ref{disc},
then $H_{\dot{\theta}}$ has purely absolutely continuous spectrum
for arbitrary $\theta$, and, in particular, $N(H_{\beta-\eps}; \E
-\lambda) = 0$ for any $\lambda > 0$.
\end{remark}

%%%%%%%%%%%%%%%%%%%%%%%%%%%%%%%%%%%%%%%%%%%%%%%%%%%%%%%%%%%%%%%%%%%%%%%%%%
\subsection{Auxiliary results}
This subsection contains some auxiliary results needed for the
proof of Theorem \ref{good}. Let $X_1$ and $X_2$ be two Hilbert
spaces. We denote by $S_{\infty}(X_1,X_2)$ the class of linear
compact  operators $T: X_1 \to X_2$. If $X_1 = X_2 = X$, we write
$S_{\infty}(X)$ instead of $S_{\infty}(X,X)$.

Let $T=T^* \in S_{\infty}(X)$. Denote by ${\mathbb P}_J(T)$ the
spectral projection of $T$ associated with  the interval $J
\subset \re$. For $s>0$ set
$$
n_{\pm}(s;T): = {\rm rank}\; {\mathbb P}_{(s,\infty)}(\pm T).
$$
If $T_j = T_j^* \in S_{\infty}(X)$, $j=1,2$, then  the Weyl
  inequalities
\begin{equation}\label{weyl}
n_{\pm}(s_1+s_2,T_1+T_2) \leq  n_{\pm}(s_1,T_1)+  n_{\pm}(s_2,T_2)
\end{equation}
hold for $s_1>0$ and $s_2>0$ (see \cite[Chapter I, Eq. (1.31)]{BSo1}).

For  $T\in S_{\infty}(X_1,X_2)$ put
\begin{equation} \label{42}
n_*(s;T): = n_+(s^2;T^*T), \quad s>0.
\end{equation}
If $T_j \in S_{\infty}(X_1,X_2)$, $j=1,2$, then  the  Ky Fan
  inequalities
\begin{equation}\label{kyfan}
n_*(s_1+s_2,T_1+T_2) \leq  n_*(s_1,T_1)+  n_*(s_2,T_2)
\end{equation}
hold for $s_1>0$ and $s_2>0$ (see \cite[Chapter I, Eq. (1.32)]{BSo1}).

Denote by $S_m(X_1,X_2)$, $m \in [1,\infty)$, the Schatten-von
Neumann classes of linear compact operators $T: X_1 \to X_2$ for
which the norm $\|T\|_m: = ({\rm Tr}\,(T^* T)^{m/2})^{1/m}$ is
finite. If $T \in S_m(X_1,X_2)$, $m \in [1,\infty)$, then the
elementary Chebyshev-type inequality
    \bel{chebyshev}
    n_*(s; T) \leq s^{-m} \|T\|_m^m
    \ee
    holds  for any $s>0$.
\begin{lemma} \label{l41} Let $J \subset \re$. Assume that $G :
{\rm L}^2(J) \to {\rm L}^2(\Omega)$ is a bounded  operator with
integral kernel $g \in {\rm L}^{\infty}(\Omega \times J)$. Let $f
\in {\rm L}^m(\Omega)$, $h \in {\rm L}^m(J)$, with $m \in
[2,\infty)$. Then the operator  $fGh \in S_m({\rm L}^2(J), {\rm
L}^2(\Omega))$, and the inequality
    \bel{smnorm}
    \|fGh\|_m^m \leq C_m  \|f\|_{{\rm L}^m(\Omega)}^m \|h\|_{{\rm L}^m(J)}^m
    \ee
holds with   $C_m =C_m(G) : = \|G\|^{m-2} \|g\|_{{\rm L}^{\infty}(\Omega \times J)}^2$.
\end{lemma}
\begin{proof}
Assume at first that $f \in
{\rm L}^{\infty}(\Omega)$, $h \in {\rm L}^{\infty}(J)$. Then, evidently, the operator $fGh$ is bounded, and we have
$$
 \|fGh\| \leq \|G\| \|f\|_{{\rm L}^{\infty}(\Omega)} \|h\|_{{\rm L}^{\infty}(J)}.
 $$
 Assume now that $f \in
{\rm L}^2(\Omega)$, $h \in {\rm L}^2(J)$. Then $fGh$ is
Hilbert-Schmidt, and we have
$$
 \|fGh\|_2 \leq
 \|g\|_{{\rm L}^{\infty}(\Omega \times J)} \|f\|_{{\rm L}^2(\Omega)} \|h\|_{{\rm L}^2(J)}.
 $$
 Applying a standard bilinear interpolation (see \cite[Section 4.4]{bl}), we get \eqref{smnorm}.
\end{proof}

\begin{remark} Results similar to Lemma \ref{l41} are contained in \cite{BSo3} and  \cite[Lemma 2.3]{R}. We include the proof of the lemma just for the convenience of the reader.
\end{remark}

Combining \eqref{chebyshev} with \eqref{smnorm}, we obtain the following

\begin{follow} \label{f41}
Under the hypotheses of Lemma \ref{l41} we have
\bel{ocenka}
    n_*(s; fGh) \leq s^{-m} C_m  \|f\|_{{\rm L}^m(\Omega)}^m \|h\|_{{\rm L}^m(J)}^m
    \ee
 for each $s>0$.
 \end{follow}

The following lemma contains standard results on the eigenvalue
asymptotics for 1D Schr\"odinger operators with decaying
attractive potentials.

\begin{lemma} \label{l42}
Assume that $V = \overline{V} \in {\rm L}^{\infty}(\re)$ satisfies
    \bel{g29}
    |V(x)| \leq C (1 + |x|)^{-\alpha}, \quad x \in \re,
    \ee
    with some constants $\alpha > 0$ and $C>0$. Let $\hbar > 0$, and
    $$
    {\mathcal H}(\hbar, V): = - \hbar^2 \frac{d^2}{dx^2} - V
    $$
    be the
    1D Schr\"odinger operator with domain ${\rm H}^2(\re)$,
    self-adjoint in ${\rm L}^2(\re)$. \\
    (i) Assume that $\alpha \in (0,2)$ and there exists a constant
$l > 0$ such that
    \bel{g28}
\lim_{|x| \to \infty} |x|^{\alpha} V(x) = l.
    \ee
Then we have
$$
\lim_{\lambda\downarrow 0}\,  \lambda^{\frac 1\alpha-\frac 12}\,
N({\mathcal H}(\hbar, V);  -\lambda) = \frac{2 l^{\frac
1\alpha}}{\pi\alpha \hbar} B\left(\frac 32,\, \frac 1\alpha-\frac
12\right).
$$
(ii) Assume that \eqref{g28} holds with $\alpha = 2$. Then we
have
$$
\lim_{\lambda\downarrow 0}\,  |\ln{\lambda}|^{-1}\, N({\mathcal
H}(\hbar, V);  -\lambda)  = \frac{1}{\pi} \left(\frac{l}{\hbar^2} -
\frac{1}{4}\right)_+^{1/2}.
    $$
If, moreover, $l < \frac{\hbar^2}{4}$, then
    \bel{g30}
N({\mathcal H}(\hbar, V);  -\lambda) = O(1), \quad \lambda
\downarrow 0.
    \ee
    (iii) Suppose that \eqref{g29} holds with
$\alpha > 2$. Then we have again \eqref{g30}.
\end{lemma}

The first part of the lemma is quite close to \cite[Theorem XIII.82]{RS4},
the proof of the second part can be found in \cite{KS1}, while the
third part follows from the result of \cite[Problem 22, Chapter XIII]{RS4}.

\subsection{Proof of Theorem \ref{good}}
The strategy of the proof of Theorem \ref{good} is to reduce the
problem of the eigenvalue asymptotics of $H_{\beta-\eps}$ to the one for an effective one-dimensional
Schr\"odinger operator
$$
- \mu \, \frac{d^2}{dx^2} - V_{eff},
$$
where  $\mu$ is  defined in Corollary \ref{gf1}, and the
effective potential  is given by
%%%%%%%%%%%%%%%%%%%%%%%%%%%%%%%%%%%%%%%%%%%%%
$$
V_{eff}(x)= 2 \beta  \|\pd_\varphi \psi_1(\cdot;0)\|_{{\rm L}^2(\omega)}^2
    \eps(x), \quad x \in \re.
    $$
    Once this is done, we use Lemma \ref{l42} to
conclude the proof.

The reduction to the one-dimensional problem is rather lengthly and
therefore we will divide it in several steps.

\subsubsection{Projection on the bottom of the essential spectrum}
Pick $\delta > 0$ so small that the eigenvalue $E_1(p)$ is simple for $p \in [-\delta, \delta]$.
As explained in at the end of Section \ref{prelim}, we assume that the mappings
$[-\delta, \delta] \ni p \mapsto E_1(p) \in (0,\infty)$ and $[-\delta, \delta] \ni p
\mapsto \psi_1(\cdot;p) \in {\rm H}^2(\omega)$ are analytic.

Introduce the orthogonal projections
$$
\pi(p) : = |\psi_1(\cdot;p)\rangle \langle\psi_1(\cdot;p)|, \quad p \in
[-\delta, \delta],
$$
acting in ${\rm L}^2(\omega)$.  Denote by $\chi_\delta : \re \to \{0,1\}$ the characteristic
function of the interval $(-\delta, \delta)$. Set
$$
{\mathcal P}_{\delta} : = \int_{\re}^\oplus \chi_\delta(p) \pi(p)
dp, \quad P = P_{\delta} : = {\mathcal F}^* {\mathcal P}_{\delta} {\mathcal
F}.
$$
Evidently, $P$  is an orthogonal projection acting in ${\rm
L}^2(\Omega)$.  Put $Q = Q_{\delta} : = I - P_{\delta}$. Since $P$
and $Q$ commute with $H_\beta^{-1/2}$, they leave ${\rm
H}_0^1(\Omega)$ invariant.

Denote by ${\mathcal Z}_1(\eps) = {\mathcal Z}_1(\eps, \beta,
\delta)$ (respectively, ${\mathcal Z}_2(\eps) = {\mathcal Z}_2(\eps,
\beta, \delta)$) the self-adjoint operator generated in the Hilbert
space $P{\rm L}^2(\Omega)$ (respectively, in $Q{\rm L}^2(\Omega)$)
by the restriction onto $P{\rm H}_0^1(\Omega)$ (respectively, onto
$Q{\rm H}_0^1(\Omega)$) of the quadratic form ${\mathcal Q}_{\beta -
\eps}$ defined in \eqref{q}. Then the mini-max principle implies
    \bel{g2}
    N({\mathcal Z}_1(\eps); \E - \lambda) \leq N(H_{\beta - \eps} ; \E -
    \lambda), \quad \lambda > 0.
    \ee
     Pick $u \in
    {\rm H}_0^1(\Omega)$ and put $v = Pu$, $w = Qu$ so that $u = v + w$.
Then we have
$$
\langle W_{\eps, \beta} u, u\rangle = : r_{\eps,\beta}[u] =
r_{\eps,\beta}[v] + r_{\eps,\beta}[w]
$$
$$
-4\beta {\rm Re} \langle G_0 v, G_0 w\rangle - 2{\rm Re} \langle
G_0 v, G_1 w\rangle - 2{\rm Re} \langle G_1 v, G_0 w\rangle +
2{\rm Re} \langle G_2 v, G_2 w\rangle
$$
where $\langle \cdot , \cdot \rangle$ denotes the scalar product
in ${\rm L}^2(\Omega)$, and $W_{\eps, \beta}$ is the operator defined in
\eqref{g0}. Next, fix $\nu \in (0,1)$, and on $P D(H_{\beta}^{1/2}) = P {\rm L}^2(\Omega)$ define the
operators
$$
G_3 : = \eps^{(1+\nu)/2}\pd_{\varphi}^2, \quad G_4 : =
\eps^{(1+\nu)/2}\pd_{\varphi}\,\pd_3, \quad G_5 : =
|\dot{\eps}|^{1/2}\pd_{\varphi}, \quad G_6 : = \eps
\pd_{\varphi}^2.
$$
Note that the operators $G_j P$, $j=3,4,5,6$, are bounded (and
compact) in ${\rm L}^2(\Omega)$. Integrating by parts, we easily
find that
$$
\langle G_0 v, G_0 w\rangle = -\langle G_3 v, \eps^{(1-\nu)/2}
w\rangle, \quad \langle G_1 v, G_0 w\rangle = -\langle G_4 v,
\eps^{(1-\nu)/2} w\rangle,
$$
$$
\langle G_0 v, G_1 w\rangle = -\langle G_4 v, \eps^{(1-\nu)/2}
w\rangle - \langle G_5 v, {\rm sign}\,\dot{\eps}
|\dot{\eps}|^{1/2} w\rangle, \quad \langle G_2 v, G_2 w\rangle = -\langle G_6 v, \eps
w\rangle.
$$
Hence, we have
 \bel{g63}
{\mathcal Q}_{\beta - \eps}[u] \geq {\mathcal Q}_{\beta - \eps}[v] -
\sum_{j=3}^6 \int_{\Omega} |G_j v|^2 d{\bf x} + {\mathcal
Q}_{\beta - \eps}[w] -
 \int_{\Omega} {\mathcal V}(x_3)|w({\bf x})|^2 d{\bf x}
    \ee
    where
    $$
    {\mathcal V}(x) : = 4(\beta^2 + 1)\eps(x)^{1-\nu} +
    |\dot{\eps}(x)| + {\eps}(x)^2, \quad x \in \re.
    $$
Denote by ${\mathcal Z}_1^+(\eps) = {\mathcal Z}_1^+(\eps, \beta)$
the self-adjoint operator generated in  $P{\rm L}^2(\Omega)$ by the
closed quadratic form
    $$
    {\mathcal Q}_{\beta - \eps}[v] -
\sum_{j=3}^6 \int_{\Omega} |G_j v|^2 d{\bf x}, \quad v \in P{\rm H}_0^1(\Omega).
$$
Similarly, denote by ${\mathcal Z}_2^+(\eps) = {\mathcal
Z}_2^+(\eps, \beta)$ the self-adjoint operator generated in  $Q{\rm
L}^2(\Omega)$ by the closed quadratic form
    $$
    {\mathcal Q}_{\beta - \eps}[w] -
 \int_{\Omega} {\mathcal V}(x_3)|w({\bf x})|^2 d{\bf x}, \quad w \in Q{\rm H}_0^1(\Omega).
$$
Then, \eqref{g63} implies
    \bel{g3}
N(H_{\beta-\eps} ; \E - \lambda) \leq
    N({\mathcal Z}_1^+(\eps) ; \E -
\lambda) + N({\mathcal Z}_2^+(\eps) ; \E -
\lambda), \quad \lambda > 0.
    \ee
Since ${\mathcal V}(x) \to 0$ as $|x| \to \infty$, and $\omega$ is a
bounded domain, we find that the multiplier by ${\mathcal V}$ is a
relative compact perturbation of $H_{\beta - \eps}$. Using this fact
and the compactness of the resolvent difference $H_{\beta
-\eps}^{-1} - H_\beta^{-1}$, we easily check that the  difference of
the resolvents of the operators ${\mathcal Z}_2^+(\eps)$ and
${\mathcal Z}_2(0)$ is a compact operator. Therefore,
$$
\inf{\sigma_{\rm ess}}({\mathcal Z}_2^+(\eps)) = \inf{\sigma_{\rm ess}}({\mathcal Z}_2(0)) =
\min{\left\{\min_{p \in \re}{E_2(p)}, \; \min_{|p| \geq \delta}
E_1(p)\right\}} > \E,
$$
and
$$
 N({\mathcal Z}_2^+(\eps) ; \E - \lambda) = O(1), \quad \lambda \downarrow 0,
 $$
which combined with \eqref{g3} implies
    \bel{g4}
    N(H_{\beta-\eps} ; \E -
\lambda) \leq N({\mathcal Z}_1^+(\eps) ; \E - \lambda) + O(1), \quad \lambda \downarrow 0.
    \ee
Fix $\eta \in (0, 2\beta)$. Recalling \eqref{heps} and \eqref{g0},
we get
    \bel{g64}
{\mathcal Q}_{\beta - \eps}[v] \leq {\mathcal Q}_{\beta}[v] - (2\beta - \eta)  \int_{\Omega} |G_0 v|^2 d{\bf x} +
\eta^{-1} \int_{\Omega} |G_1 v|^2 d{\bf x} + \int_{\Omega} |G_2 v|^2 d{\bf x}, \quad v \in P{\rm H}_0^1(\Omega),
    \ee
    $$
    {\mathcal Q}_{\beta - \eps}[v] - \sum_{j=3}^6 \int_{\Omega} |G_j v|^2 d{\bf x} \geq
    $$
    \bel{g65}
    {\mathcal Q}_{\beta}[v] - (2\beta + \eta) \int_{\Omega} |G_0 v|^2 d{\bf x}
    - \eta^{-1} \int_{\Omega} |G_1 v|^2 d{\bf x}
    - \sum_{j=3}^6 \int_{\Omega} |G_j v|^2 d{\bf x}, \quad v \in P{\rm H}_0^1(\Omega).
    \ee
    Denote by $\tilde{{\mathcal Z}}_1^-(\eps) = \tilde{{\mathcal Z}}_1^-(\eps, \beta, \eta)$
    the self-adjoint operator generated in
     $P{\rm L}^2(\Omega)$ by
    the closed quadratic form
    $$
    {\mathcal Q}_{\beta}[v] - (2\beta - \eta) \int_{\Omega} |G_0 v|^2 d{\bf x} +
    \eta^{-1} \int_{\Omega} |G_1 v|^2d{\bf x} + \int_{\Omega} |G_2 v|^2 d{\bf x}, \quad v \in P{\rm H}_0^1(\Omega).
$$
    Similarly, denote by $\tilde{{\mathcal Z}}_1^+(\eps) = \tilde{{\mathcal Z}}_1^+(\eps, \beta, \eta)$ the
    self-adjoint operator generated in  $P{\rm L}_2(\Omega)$ by
    the closed quadratic form
    $$
    {\mathcal Q}_{\beta}[v] -(2\beta + \eta) \int_{\Omega} |G_0 v|^2 d{\bf x}
    - \eta^{-1} \int_{\Omega} |G_1 v|^2 d{\bf x} -
\sum_{j=3}^6 \int_{\Omega} |G_j v|^2 d{\bf x}, \quad v \in P{\rm H}_0^1(\Omega).
$$
Then \eqref{g64} - \eqref{g65} implies
    \bel{g66}
    {\mathcal Z}_1^+ \geq \tilde{{\mathcal Z}}_1^+, \quad {\mathcal Z}_1 \leq \tilde{{\mathcal Z}}_1^-.
    \ee
    Now estimates \eqref{g2}, \eqref{g4}, and \eqref{g66} entail
\bel{g6}
N(\tilde{{\mathcal Z}}_1^-- \E + \lambda ; 0) \leq
N(H_{\beta - \eps}; \E - \lambda) \leq
    N(\tilde{{\mathcal Z}}_1^+ - \E + \lambda ; 0) + O(1),
\quad \lambda \downarrow 0.
    \ee
    The last equation shows that the eigenvalue asymptotics
of $H_{\beta - \eps}$ is determined by the asymptotics of the
reduced operator $PH_{\beta - \eps}P $ modulo some error terms.

\subsubsection{Reduction to a one-dimensional problem}
We introduce the operator $U: {\rm L}^2(-\delta, \delta) \to P{\rm
L}^2(\Omega)$ which acts on $f \in  {\rm L}^2(-\delta, \delta)$ as
follows
\begin{equation} \label{U}
(U f)(x_t,x_3):= {\mathcal F}^* \tilde{f}, \quad \tilde{f}(x_t,p) : =
%\left\{
\Bigg\{
\begin{array} {l}
\psi_1(x_t,p) f(p) \quad {\rm if} \quad x_t \in \omega, \; p \in (-\delta, \delta),\\
 \\
0 \quad {\rm if} \quad x_t \in \omega, \; p \in \re\setminus(-\delta, \delta).
\end{array}
%\right.
\end{equation}
 Then $U$ is a unitary operator and using \eqref{U} it
can be directly verified that
$$
\tilde{{\mathcal Z}}_1^-- \E + \lambda = U (M(\lambda) - (2\beta - \eta)\Gamma_0^* \Gamma_0 + \eta^{-1}
\Gamma_1^*\Gamma_1 + \Gamma_2^* \Gamma_2) U^*,
$$
$$
\tilde{{\mathcal Z}}_1^+- \E + \lambda =
U (M(\lambda) - (2\beta + \eta)\Gamma_0^* \Gamma_0 - \eta^{-1}
\Gamma_1^*\Gamma_1 - \sum_{j=3}^6 \Gamma_j^* \Gamma_j) U^*,
$$
where $M(\lambda)$ is the multiplier by $E_1(p) - \E + \lambda$ in
${\rm L}^2(-\delta, \delta)$, and
 $\Gamma_j : {\rm L}^2(-\delta, \delta) \to {\rm L}^2(\Omega)$ are
integral operators with kernels $(2\pi)^{-1/2} e^{ix_3 p}
\gamma_j(\bx,p)$, $\bx = (x_t, x_3) \in \Omega$, $p \in (-\delta,
\delta)$, defined by
$$
\gamma_0(\bx,p) : =  \eps(x_3)^{1/2} \pd_{\varphi}\psi_1(x_t;p),
\quad
 \gamma_1(\bx,p) : = i \eps(x_3)^{1/2} \psi_1(x_t;p) p,
 $$
 $$
 \gamma_2(\bx,p) : =  \eps(x_3) \pd_{\varphi}\psi_1(x_t;p), \quad
\gamma_3(\bx,p) : =  \eps(x_3)^{(1+\nu)/2} \pd_{\varphi}^2
\psi_1(x_t;p),
$$
$$
\quad \gamma_4(\bx,p) : = i \eps(x_3)^{(1+\nu)/2}
\pd_{\varphi}\psi_1(x_t;p)p, \quad \gamma_5(\bx,p) : =
|\dot{\eps}(x_3)|^{1/2} \pd_{\varphi}\psi_1(x_t;p),
$$
$$
 \gamma_6(\bx,p) : =  \eps(x_3) \pd_{\varphi}^2
\psi_1(x_t;p).
$$
Then \eqref{g6} implies
$$
N(M(\lambda) - (2\beta - \eta)\Gamma_0^* \Gamma_0 + \eta^{-1}
\Gamma_1^*\Gamma_1 + \Gamma_2^* \Gamma_2; 0) \leq
$$
$$
N(H_{\beta-\eps} ; \E - \lambda) \leq
    $$
    \bel{g6a}
N(M(\lambda) - (2\beta + \eta)\Gamma_0^* \Gamma_0 - \eta^{-1}
\Gamma_1^*\Gamma_1 - \sum_{j=3}^6 \Gamma_j^* \Gamma_j ; 0) +
O(1), \quad \lambda \downarrow 0.
    \ee
    Set
    $$
    a_{\lambda}(p) : = (E_1(p) - \E + \lambda)^{-1/2}, \quad \lambda >
    0, \quad p \in [-\delta, \delta].
    $$
    Applying the Birman-Schwinger principle, the Weyl
    inequalities \eqref{weyl}, and definition \eqref{42}, we find that for each $s \in (0,1)$ we have
    $$
N(M(\lambda) - (2\beta - \eta)\Gamma_0^* \Gamma_0 + \eta^{-1}
\Gamma_1^*\Gamma_1 + \Gamma_2^* \Gamma_2; 0) =
$$
$$
n_+(1; a_{\lambda}((2\beta - \eta)\Gamma_0^* \Gamma_0 - \eta^{-1}
\Gamma_1^*\Gamma_1 - \Gamma_2^* \Gamma_2) a_{\lambda}) \geq
$$
    \bel{g31}
    n_*(\sqrt{(1+s)/(2\beta - \eta)}; \Gamma_0 a_{\lambda}) -
n_*(\sqrt{\eta s/2}; \Gamma_1 a_{\lambda}) - n_*(\sqrt{s/2};
\Gamma_2 a_{\lambda}),
    \ee
    $$
N(M(\lambda) - (2\beta  + \eta)\Gamma_0^* \Gamma_0 - \eta^{-1}
\Gamma_1^*\Gamma_1 - \sum_{j=3}^6\Gamma_j^* \Gamma_j; 0) =
$$
$$
n_+(1; a_{\lambda}((2\beta + \eta)\Gamma_0^* \Gamma_0 + \eta^{-1}
\Gamma_1^*\Gamma_1 + \sum_{j=3}^6 \Gamma_1^* \Gamma_1)
a_{\lambda}) \leq
$$
    \bel{g32}
    n_*(\sqrt{(1-s)/(2\beta + \eta)}; \Gamma_0 a_{\lambda}) +
n_*(\sqrt{\eta s/5}; \Gamma_1 a_{\lambda}) + \sum_{j=3}^6
n_*(\sqrt{s/5}; \Gamma_j a_{\lambda}).
    \ee
    Note that on the right-hand sides of \eqref{g31} and \eqref{g32}
    there are just linear combinations of terms of the form $n_*(r; \Gamma_j
    a_{\lambda})$ with $r>0$ independent of $\lambda$ and
    $j=0,\ldots, 6$. The rest of the proof of Theorem \ref{good}
    reduces to the  asymptotic analysis as $\lambda \downarrow 0$ of these terms.
    Our aim is to show that only the ones corresponding to
the operator $\Gamma_0$ contribute to the main asymptotic term of
$N(H_{\beta-\eps};\E - \lambda)$ as $\lambda \downarrow 0$.
    First, we show that
    \bel{g10}
    n_*(r; \Gamma_j a_{\lambda}) = O(1), \quad \lambda \downarrow 0,
    \quad j=1,4,
    \ee
for every $r>0$. To this end, it suffices to apply \eqref{ocenka}
with $J = (-\delta, \delta)$, $f(\bx) = \eps(x_3)^{1/2}$ if $j=1$,
$f(\bx) = \eps(x_3)^{(1+\nu)/2}$ if $j=4$, $g(\bx, p) =
(2\pi)^{-1/2} e^{ix_3p} \pd_{\varphi}\psi_1(x_t; p)$, $h(p) = ip
a_{\lambda}(p)$, and $m \in [2,\infty)$ large enough.

Further, for $j=0,2,3,5,6$ we define the operators
$\tilde{\Gamma}_j : {\rm L}^2(\re) \to {\rm
L}^2(\Omega)$ as integral operators whose kernels are obtained by
substituting $\psi_1(x_t; 0)$ for $\psi_1(x_t; p)$ in the
expressions for the integral kernels of $\Gamma_j$.
Denote by $\tilde{\Gamma}_{j,\delta}$ the restriction of $\tilde{\Gamma}_j$ onto $L^2(-\delta, \delta)$.
Set
$$
\tilde{a}_{\lambda}(p) = (\mu p^2 + \lambda)^{-1/2}, \quad p \in
\re, \quad \lambda > 0.
$$
Pick $s \in (0,1)$, and bearing in mind \eqref{nd2} choose $\delta > 0$ so small that we have
    \bel{g33}
    (1+s)^{-1} \tilde{a}_{\lambda}(p) \leq a_\lambda(p) \leq (1-s)^{-1}
    \tilde{a}_{\lambda}(p), \quad p \in [-\delta, \delta], \quad
    \lambda > 0.
    \ee
    Our next goal is to show that the quantities $n_*(r; \Gamma_j
    a_\lambda)$, $r>0$, $j=2,3,5,6$, appearing at the right-hand
    sides of \eqref{g31} and \eqref{g32} are bounded under the hypotheses of Theorem \ref{good} (ii) - (iii),
    and do not contribute to the
    main asymptotic term as $\lambda \downarrow 0$ of $N(H_{\beta-\eps};
    \E - \lambda)$ under the hypotheses of Theorem \ref{good} (i), even though in the last case
    they might not be bounded in contrast to the cases $j=1,4$. The upper bound  in
    \eqref{g33} combined with the mini-max principle, and the Ky Fan
    inequalities \eqref{kyfan} imply
    \bel{g34}
    n_*(r; \Gamma_j a_\lambda) \leq n_*(r(1-s)^2; \tilde{\Gamma}_j
    \tilde{a}_\lambda) + n_*(rs; (\Gamma_j -  \tilde{\Gamma}_{j, \delta})
    a_\lambda), \quad r>0, \quad j=2,3,5,6.
    \ee
    It is quite easy to see that
    \bel{g35}
n_*(rs; (\Gamma_j - \tilde{\Gamma}_{j, \delta})
    a_\lambda) = O(1), \quad \lambda \downarrow 0,
    \ee
    for $r>0$ and $j=2,3,5,6$; for example, if $j=2$ it suffices to apply
    \eqref{ocenka} with $f(\bx) = \eps(x_3)$, $g(\bx, p) =
    (2\pi)^{-1/2} e^{ix_3p} p^{-1} (\pd_\varphi\psi_1(x_t;p) -
    \pd_\varphi\psi_1(x_t;0))$, $h(p) = a_\lambda(p) p$, and $m \in
    [2, \infty)$ large enough, and if $j=3,5,6$, the argument is the
    same with appropriate choice of $f$, $g$, and $h$.

    Now, the Birman-Schwinger principle implies that for each
    $r>0$ we have
    \bel{g36}
    n_*(r; \tilde{\Gamma}_j \tilde{a}_\lambda) =
    n_+(r^2; \tilde{a}_\lambda \tilde{\Gamma}_j^* \tilde{\Gamma}_j
    \tilde{a}_\lambda) = N({\mathcal H}(\sqrt{\mu}, r^{-2} V_j);
    -\lambda), \quad j=2,3,5,6, \quad \lambda>0,
    \ee
    where ${\mathcal H}$ is the 1D Schr\"odinger operator defined
    in Lemma \ref{l42}, and
    $$
    V_2(x) : =  \|\pd_\varphi \psi_1(\cdot;0)\|_{{\rm L}^2(\omega)}^2 \eps(x)^2, \quad
    V_3(x) : =  \|\pd^2_\varphi \psi_1(\cdot;0)\|_{{\rm L}^2(\omega)}^2 \eps(x)^{1+\nu},
    $$
    $$
    V_5(x) : =  \|\pd_\varphi \psi_1(\cdot;0)\|_{{\rm L}^2(\omega)}^2 |\dot{\eps}(x)|,
    \quad V_6(x) : =  \|\pd^2_\varphi \psi_1(\cdot;0)\|_{{\rm L}^2(\omega)}^2 \eps(x)^2, \quad x \in \re.
    $$
    Assumption \ref{estimate} implies that $\lim_{|x|\to\infty}|x|^{\alpha}
    V_j(x) = 0$, $j = 2,3,5,6$. Applying Lemma \ref{l42} to the counting function
    $N({\mathcal H}(\sqrt{\mu}, r^{-2} V_j); -\lambda)$, and putting together \eqref{g34} -- \eqref{g36}, we
    obtain that for $j=2,3,5,6$, and $r>0$, we have
    \bel{g37}
n_*(r; \Gamma_j a_\lambda) = \left\{
    \begin{array} {l}
o(\lambda^{\frac{1}{2} - \frac{1}{\alpha}}) \; \text{if Assumption
\ref{asymptotics} holds with} \; \alpha \in
    (0,2),\\
     O(1) = o(|\ln{\lambda}|) \; \text{if Assumption \ref{asymptotics} holds with} \; \alpha = 2,\\
     O(1) \; \text{if Assumption \ref{estimate} holds with} \; \alpha >
     2,
    \end{array}
 \right.
    \ee
as $\lambda \downarrow 0$.
\subsubsection{Eigenvalue asymptotics for the one-dimensional operator}
Finally, we will show that the quantities $n_*(r; \Gamma_0
    a_\lambda)$, $r>0$,   on  the right-hand
    sides of \eqref{g31} and \eqref{g32}, generate the
    main asymptotic term as $\lambda \downarrow 0$ of $N(H_{\beta-\eps};
    \E - \lambda)$ under the hypotheses of Theorem \ref{good} (i) - (ii),
    or are bounded under the hypotheses of Theorem \ref{good}
    (iii).
    Similarly to \eqref{g34}, the estimates  \eqref{g33}, combined with the
    mini-max principle, and the  Ky Fan
    inequalities \eqref{kyfan} imply
    $$
n_*(r(1+s)^3; \tilde{\Gamma}_0
    \tilde{a}_\lambda) - n_*(rs(1+s)^2;\tilde{\Gamma}_0
    \tilde{a}_\lambda w_\delta) - n_*(rs; (\Gamma_0 - \tilde{\Gamma}_{0,\delta})
    a_\lambda)\leq
    $$
    \bel{g38}
    n_*(r; \Gamma_0 a_\lambda) \leq n_*(r(1-s)^2; \tilde{\Gamma}_0
    \tilde{a}_\lambda) + n_*(rs; (\Gamma_0 - \tilde{\Gamma}_{0,\delta})
    a_\lambda), \quad r>0,
    \ee
    where $w_\delta : \re \to \{0,1\}$ denotes the characteristic function of the set $\re \setminus(-\delta,
    \delta)$.

    By analogy with \eqref{g35}, we find that
    \bel{g39}
n_+(r; (\Gamma_0 - \tilde{\Gamma}_{0,\delta})
    a_\lambda) = O(1), \quad \lambda \downarrow 0, \quad r>0.
    \ee
    Further, applying \eqref{ocenka} with $$
    J = \re, \quad  f(\bx) = \eps(x_3)^{1/2}, \quad g(\bx, p) =
    (2\pi)^{-1/2} e^{ix_3p} \pd_\varphi\psi_1(x_t;0),  \quad h(p) = a_\lambda(p) w_\delta(p),
    $$
     and $m \in
    [2, \infty)$ large enough, we get
    \bel{g40}
n_*(r;\tilde{\Gamma}_0
    \tilde{a}_\lambda w_\delta)
     = O(1), \quad \lambda \downarrow 0, \quad r>0.
    \ee
    The Birman-Schwinger principle implies that for each
    $r>0$ we have
    \bel{g41}
    n_*(r; \tilde{\Gamma}_0 \tilde{a}_\lambda) =
    n_+(r^2; \tilde{a}_\lambda \tilde{\Gamma}_0^* \tilde{\Gamma}_0
    \tilde{a}_\lambda) = N({\mathcal H}(\sqrt{\mu}, r^{-2} V_0);
    -\lambda), \quad \lambda>0,
    \ee
    where ${\mathcal H}$ is the 1D Schr\"odinger operator defined
    in Lemma \ref{l42}, and
    $$
    V_0(x) : =  \|\pd_\varphi \psi_1(\cdot;0)\|_{{\rm L}^2(\omega)}^2
    \eps(x).
    $$
    Applying Lemma \ref{l42}, and bearing in mind \eqref{g38} -- \eqref{g40}, we find
    that
    \bel{g42}
\lim_{\lambda\downarrow 0}\,  \lambda^{\frac 1\alpha-\frac 12}\,
n_*(r; \Gamma_0 a_\lambda) = \frac{2 (r^{-2} \|\pd_\varphi
\psi_1(\cdot;0)\|_{{\rm L}^2(\omega)}^2 L)^{\frac
1\alpha}}{\pi\alpha \sqrt{\mu}} B\left(\frac 32,\, \frac
1\alpha-\frac 12\right).
    \ee
if Assumption \ref{asymptotics} with $\alpha \in (0,2)$ holds
true,
    \bel{g44}
\lim_{\lambda\downarrow 0}\,  |\ln{\lambda}|^{-1} n_*(r; \Gamma_0
a_\lambda) = \frac{1}{\pi} \left(\frac{r^{-2} \|\pd_\varphi
\psi_1(\cdot;0)\|_{{\rm L}^2(\omega)}^2 L}{\mu} -
\frac{1}{4}\right)_+^{1/2},
    \ee
    if Assumption \ref{asymptotics} with $\alpha = 2$ holds
true, and
    \bel{g45}
    n_*(r; \Gamma_0 a_\lambda) = O(1), \quad \lambda \downarrow 0,
    \ee
    if Assumption \ref{estimate} holds true with $\alpha > 2$.
If, moreover, Assumption \ref{asymptotics} holds with $\alpha =
2$, and  $r^{-2} \|\pd_\varphi \psi_1(\cdot;0)\|_{{\rm
L}^2(\omega)}^2 L < \frac{\mu}{4}$, then we have \eqref{g45}.

Since the numbers $s \in (0,1)$ and $\eta \in (0,2\beta)$ in
\eqref{g31} -- \eqref{g32} could be chosen arbitrarily small, we
find that \eqref{g6a} -- \eqref{g32}, \eqref{g37}, \eqref{g38},
\eqref{g42} -- \eqref{g45} imply that under the appropriate
assumptions of Theorem \ref{good}, asymptotic relations
\eqref{g26}, \eqref{g27}, or \eqref{g1} hold true.\\

{\bf Acknowledgements}. Hynek Kova\v r\'{\i}k was partially supported by the German Research
Foundation (DFG) under Grant KO 3636/1-1. Georgi Raikov was partially supported by the Chilean
Scientific Foundation {\em Fondecyt} under Grant 1050716.

\bigskip

{\sc Ph. Briet}\\
Centre de Physique Th\'eorique\\
CNRS-Luminy, Case 907\\
13288 Marseille, France\\
E-mail: briet@cpt.univ-mrs.fr\\

{\sc H. Kova\v r\'{\i}k}\\
Universit\`a degli Studi di Modena e Reggio Emilia\\
Dipartimento di Matematica\\
 Via Campi 213/B\\
I-41100 Modena, Italy\\
E-mail: hynek.kovarik@unimore.it \\

{\sc G. Raikov}\\
%Departamento de Matem\'aticas\\
Facultad de Matem\'aticas\\
Pontificia Universidad Cat\'olica de Chile\\
Av. Vicu\~na Mackenna 4860\\ Santiago de Chile\\
E-mail: graikov@mat.puc.cl\\

{\sc E. Soccorsi}\\
Centre de Physique Th\'eorique\\
CNRS-Luminy, Case 907\\
13288 Marseille, France\\
E-mail: soccorsi@cpt.univ-mrs.fr\\

\end{document}